\documentclass[a4paper, 11pt]{article}

\usepackage{amssymb}
\usepackage[definethebibliography]{easybib}
\usepackage{amsthm}
\usepackage{amsmath}
\usepackage{mathrsfs}

\begin{document}

\newtheorem{thm}{Theorem}[section]
\newtheorem{lem}[thm]{Lemma}
\newtheorem{rem}[thm]{Remark}
\newtheorem{cor}[thm]{Corollary}
\newtheorem{prop}[thm]{Proposition}


\renewcommand{\theequation}{\arabic{section}.\arabic{equation}}

\def\proof{\noindent{\it Proof.\ }}
\newcommand{\RR}{\mathbb R}
\newcommand{\be}{\begin{equation}} 
\newcommand{\ee}{\end{equation}}
\newcommand{\bea}{\begin{eqnarray}} 
\newcommand{\eea}{\end{eqnarray}}
\newcommand{\bean}{\begin{eqnarray*}} 
\newcommand{\eean}{\end{eqnarray*}}
\newcommand{\rf}[1]{(\ref {#1})}
\newcommand{\un}{{\rm 1\!I}}
\newcommand{\dx}{\,{\rm d}x}
\newcommand{\dy}{\,{\rm d}y}
\newcommand{\ds}{\,{\rm d}s}
\newcommand{\dt}{\,{\rm d}t}
\newcommand{\la}{\langle}
\newcommand{\ra}{\rangle} 
\def\t{\tau}
\def\s{\sigma} 
\def\e{\varepsilon}
\def\p{\partial}
\def\div{\nabla\cdot}
\def\f{\varphi}
\def\r{\varrho}
\def\h{{H^1(\Omega)}}

\newcommand{\question}[1]{{\par\medskip\hrule\medskip\noindent {\small\sf #1} 
\medskip\hrule\medskip\par}}
\def\X{\mathcal{X}}

\title{\Large\textbf{Blow up of solutions \\ to  generalized Keller--Segel model}} 
\author{Piotr Biler \qquad  Grzegorz Karch\\
\small   Instytut Matematyczny,
Uniwersytet Wroc{\l}awski \\
\small pl. Grunwaldzki 2/4, 50--384 Wroc{\l}aw, Poland\\
\small{\tt \{Piotr.Biler,Grzegorz.Karch\}@math.uni.wroc.pl}}

\date{\today}
\maketitle

\begin{abstract} 
\noindent The existence and 
nonexistence of global in time solutions is studied for a class of 
equations generalizing the chemotaxis model of Keller and Segel. These equations 
involve L\'evy diffusion operators and general potential type nonlinear 
terms. 
\end{abstract}

\noindent{\bf Key words and phrases:} nonlocal parabolic equations, blow up of 
solutions, L\'evy diffusion, chemotaxis, moment method
\bigskip

\noindent {\sl 2000 Mathematics Subject Classification:} 35Q, 35K55, 35B40. 
\baselineskip=18.8pt

\section{Introduction} 

\setcounter{equation}{0}

We consider in this paper the following nonlinear nonlocal  evolution 
equation generalizing the well known Keller--Segel model of chemotaxis 
\be
\partial_tu+(-\Delta)^{\alpha/2}u+\nabla\cdot(uB(u))=0, \label{eq} 
\ee
for $(x,t)\in \RR^d\times\RR^+$, where the anomalous diffusion is 
modeled by a fractional power of the Laplacian, $\alpha\in(1,2)$, and the linear (vector) operator $B$ is defined (formally) as 
\be 
B(u)=\nabla((-\Delta)^{-\beta/2}u).\label{B}
\ee
For  $\beta\in(1,d]$, and $d\ge 2$, one may express the nonlocal nonlinearity in \rf{eq} using convolution operators since 
\be 
B(u)(x)=
s_{d,\beta}\int_{\RR^d}\frac{x-y}{|x-y|^{d-\beta+2}}u(y) \dy,\label{BB}
\ee 
with some $s_{d,\beta}>0$, and the assumption $\beta>1$ is needed for the convergence of this integral. 

Of course, the choice $\alpha=2$, $\beta=2$ in \rf{eq}--\rf{BB} corresponds
 to the usual Keller--Segel system studied mainly in space dimensions $d=1,\, 2,\, 3$. 
It is well known that if $\beta=2$, the one-dimensional system \rf{eq}--\rf{BB} 
possesses global in time solutions not only in the case of classical Brownian 
diffusion $\alpha=2$ but also in the fractional diffusion case $1<\alpha<2$ as 
was shown by C.~Escudero in \cite{E}. 
On the other hand, there are many results on the nonexistence of global in time 
solutions with ``large'' initial data if $d\ge2$ and $\alpha=2$, $\beta=2$, see, 
{\it e.g.}, \cite{JL,B3,B-AMSA,BDP}. Even if $d\ge 2$, $\alpha=2$ and $1<\beta\le d$, 
there are results on the blow up of solutions with suitably chosen initial data, 
see \cite{BW} (caution: the notation in \cite[Prop. 4.2]{BW} differs from that 
in the present paper).  
Let us finally recall that, in the limit case $\alpha=2$, $\beta=d$, mass 
$M=\int_{\RR^d}u_0(x)\dx$ is the critical parameter for the blow up, see 
\cite[Prop. 4.1]{BW} and \cite{CPS}. 

The usual method of proving the nonexistence of global in time nonnegative
and nontrivial solutions,
used in the abovementioned papers, consists in the study the evolution of
the second moment of a solution $w_2(t)=\int_{\RR^d}|x|^2u(x,t)\dx$ and to show
(via suitable differential inequalities) that $w_2(t)$ vanishes 
for some $t>0$.
The second moment of a typical solution to an evolution
 equation with  fractional Laplacian  cannot be finite, see {\it  e.g.} \cite{BK08}.
Hence, our goal in this paper is to generalize the classical virial method 
and
to show  the blow up of solution  system \rf{eq}--\rf{BB} by studying moments of lower order $\gamma \in(1,2)$
\be
w_\gamma=\int_{\RR^d} |x|^\gamma u(x)\dx.\label{m-gamma}
\ee

After this paper was completed, we discovered a recent preprint \cite{LRZ08}
where the authors show the blow up of solution to system \rf{eq}--\rf{BB} with fractional diffusion in the particular case $d=2$ and $\beta=2$. Our argument  is different than that in  \cite{LRZ08}, shorter, seems to be more direct, and applies in more general situations. Moreover, we are able to formulate a simple condition on the initial data which leads to the blow up
in a finite time  of the corresponding solution. 

\medskip

{\bf Notation.}
The $L^p$-norm   of a Lebesgue
measurable, real-valued function $v$ defined on $\RR^d$
is denoted by $\|v\|_p$.
The constants (always independent of $x,t$) will be
denoted by the same letter $C$, even if they may vary from line to line.
Sometimes, we write, {\it e.g.},  $C=C(*)$ when we want to
emphasize the dependence of $C$ on a~parameter ``$*$''.

\section{Main results}
\setcounter{equation}{0}

The crucial role in the approach 
in this paper is played by the following scaling property of system \rf{eq}--\rf{BB}
\begin{equation}\label{scal}
u^\lambda(x,t) =\lambda^{\alpha+\beta-2}u(\lambda x, \lambda^\alpha t)
\quad \mbox{for all}\quad \lambda>0, 
\end{equation}
in the sense that if $u$ is a solution to \rf{eq}--\rf{BB}, then $u^\lambda$ 
is so. In particular, in our construction of solutions to  \rf{eq}--\rf{BB}  
we use the fact that
the usual norm of the Lebesgue space $L^{d/(\alpha+\beta-2)}(\RR^d)$ 
is invariant under the transformation 
$u_0(x)\mapsto \lambda^{\alpha+\beta-2}u_0(\lambda x)$ for every $\lambda>0$.

\begin{thm}\label{existence} 
Assume that $d\ge 2$, $\alpha\in(1,2]$, and  $\beta\in (1,d]$.
Let 
$$ \max\left\{\frac{d}{\alpha+\beta-2}, \frac{2d}{d+\beta-1} \right\}<p\leq  d.
$$
\begin{itemize}
\item[i)]
For every $u_0\in L^p(\RR^d)$ there exists $T=T(\|u_0\|_p)$ and the unique local in time mild solution $u\in C([0,T], L^p(\RR^d))$ of system \rf{eq}--\rf{BB} with $u_0$ as the initial condition. 

\item[ii)] 
There is $\varepsilon>0$ such that for every $u_0\in L^{d/(\alpha+\beta-2)}(\RR^d)$ satisfying 
\be
\|u_0\|_{d/(\alpha+\beta-2)}\leq\varepsilon\label{smallness}
\ee 
there exists a global in time mild solution $u\in C([0,\infty), L^p(\RR^d))$ of system \rf{eq}--\rf{BB} with $u_0$ as the initial condition.
\end{itemize}
Moreover, if $u_0(x)\ge 0$, then the solution $u$ in  either i) or ii) above 
is nonnegative. Finally, if  $u_0\in L^1(\RR^d)$, then the corresponding 
solution   
conserves mass 
\be
\int_{\RR^d}u(x,t)\dx= \int_{\RR^d}u_0(x)\dx\equiv M.\label{mass}
\ee 
\end{thm}

\medskip

Recall that $\alpha>1$ is a usual assumption 
(\cite[Th. 2.2]{BW}, \cite{BKW1,BKW2}) 
which permits us to control locally the nonlinearity in \rf{eq}--\rf{BB} 
by the linear term.

The results stated above 
 can be easily generalized for equations with general L\'evy 
diffusion operators considered in \cite{BKW1,BKW2} but we do not pursue this 
question here. 
We refer the reader to the above mentioned papers, as well as \cite{BW}, for 
physical motivations to study such equations. On the other hand, motivations 
stemming from probability theory ({\sl propagation of chaos} property for 
interacting particle systems) can be found in, {\it e.g.}, \cite{BFW}. 

\begin{rem}
{\rm For $\alpha\in(1,2)$ and $\beta>1$ satisfying $\alpha+\beta>d+2$, 
it can be shown that local in time solutions 
can be continued to the global in time ones, see \cite[Th.~3.2]{BW}. 
However, in \cite{BW}, another approach (via weak solutions) 
is used to construct solutions of system \rf{eq}--\rf{BB}. }
\end{rem}

The proof of Theorem \ref{existence} on local and global solutions 
to \rf{eq}--\rf{BB} follows a more or less standard reasoning which we sketch 
in Section \ref{sec:exist}.
Our main goal, however, is to prove the finite time blow up of solutions to 
the nonlocal system \rf{eq}--\rf{BB}.

\begin{thm}\label{main}
Assume that $d\geq 2$. The solution of \rf{eq}--\rf{BB} with 
a~non\-negative and nonzero initial condition 
$u(x,0)=u_0(x)$ 
blows up in a finite time in each of the following cases:

\begin{itemize}
\item[i)] {\em {(large mass)}} for $\alpha=2$, $\beta=d$, 
$u_0\in L^1(\RR^d, (1+|x|^2)\dx)$,
and if 
 $$M=\int_{\RR^d}u_0(x)\dx>2d/s_{d,\beta},$$
with the constant $s_{d,\beta}$ defined in  \rf{BB}; 
in particular: $s_{2,2}=\frac{1}{2\pi}$ so that the threshold value of $M$ is $8\pi$ if $d=2$;

\item[ii)] {\em {(high concentration)}} 
for $\alpha\in(1,2]$ and  $\beta\in(1,d]$ satisfying $\alpha+\beta<d+2$,
$u_0\in L^1(\RR^d, (1+|x|^\gamma)\dx)$ for some $\gamma\in (1,\alpha)$,
and if
\begin{equation}
\frac{\int_{\RR^d} |x|^\gamma u_0(x) \dx}{\int_{\RR^d}u_0(x)\dx}
\leq c \Bigg(\int_{\RR^d}u_0(x)\dx\Bigg)^{\frac{\gamma}{d+2-\alpha-\beta}}\label{cond:blow-up}
\end{equation}
for certain (sufficiently small) constant $c>0$ independent of $u_0$.
\end{itemize}
\end{thm}

The result stated in i) is essentially contained in \cite{CPS}.  
The condition for blow up in the form \rf{cond:blow-up} appeared already in \cite{B3} and \cite{N2000}, 
of course, for $\alpha=2$ only. 
Note that  i) is a limit case of ii). Indeed, \rf{cond:blow-up} written as 
$$
\left(\frac{\int_{\RR^d} |x|^\gamma u_0(x) \dx}{\int_{\RR^d}u_0(x)\dx}\right)^{\frac{d+2-\alpha-\beta}{\gamma}}\le cM,
$$ 
becomes a condition on (sufficiently large) mass: $1\le cM$, when $(\alpha+\beta)\nearrow (d+2)$.

We have to emphasize that Theorem \ref{main}.ii contains a result  which is 
new even for the classical parabolic-elliptic Keller--Segel model 
({\it i.e.}~equations~\rf{eq}--\rf{BB} 
with $\alpha=\beta =2$). Indeed, for $d\geq 3$,
the conditions from part ii) of Theorem \ref{main} guarantee the blow
up in a finite time if the moment of order $\gamma$ of the initial condition  is 
finite for some $\gamma \in (1,2)$.  All other known proofs of
the blow up required just $\gamma=2$. An analogous result for $d=2$ is stated below.

\begin{cor}\label{cor:main}
Assume that $\alpha=\beta=d=2$ in equations \rf{eq}--\rf{BB}. Suppose that
there exists $\gamma\in (1,2)$ such that 
$u_0 \in  L^1(\RR^2, (1+|x|^\gamma)\dx)$.
There exists $M_\gamma>0$ 
such that if 
$$
\int_{\RR^2}u_0(x)\dx> M_\gamma, \quad
$$
then the solution of \rf{eq}--\rf{BB}  with the initial condition 
$u(x,0)=u_0(x)$ blows up in a finite time. 
\end{cor}

It is conjectured that $M_2=M_\gamma=8\pi$ for every $\gamma\in (1,2)$ and we
expect that it can be shown by  a careful analysis of constants appearing 
in inequalities \rf{der:est}--\rf{grad} and in \rf{phi-phi}--\rf{evo} in the proof of Theorem  \ref{main}.
Here, one should recall that {\em radially symmetric} solutions to
 equations \rf{eq}--\rf{BB} with $\alpha=\beta=d=2$ blow up in finite time
under the assumption that their mass is larger than $8\pi$ and no moment 
condition imposed on the initial data is necessary, 
see \cite[Prop. 2.2, Th. 3.1]{B-BCP} for a detailed presentation. 

\medskip

\begin{rem}
{\rm Due to the translation invariance of problem \rf{eq}--\rf{BB}, the conditions on moments can be imposed on the quantity 
$$\inf_{x_0\in\RR^d}\int_{\RR^d} |x-x_0|^\gamma u_0(x)\dx=\int_{\RR^d} |x-\bar x|^\gamma u_0(x)\dx,$$
where $\bar x=(\int_{\RR^d} xu_0(x)\dx)/(\int_{\RR^d}u_0(x)\dx)$ is the center of mass of $u_0$. }
\end{rem}

\begin{rem}
{\rm 
Let us observe that the assumptions \rf{smallness} and \rf{cond:blow-up} are in a~sense complementary due to the following  elementary inequality involving 
the $L^p$-norms, mass $M=\int_{\RR^d}u(x)\dx$, and the moment  
$w_\gamma=\int_{\RR^d} |x|^\gamma u(x)\dx$ of a nonnegative function $u$ 
\be
\|u\|_p\ge CM\left(\frac{M}{w_\gamma}\right)^{\frac{d}{\gamma}\left(1-\frac{1}{p}\right)}. \label{m-L}
\ee 
To prove \rf{m-L} observe that $w_\gamma\ge R^\gamma\int_{\RR^d\setminus B_R(0)}u(x)\dx$, so that 
$$\int_{B_R(0)}u(x)\dx=M-\int_{\RR^d\setminus B_R(0)}u(x)\dx\ge M-R^{-\gamma}w_\gamma.$$ 
Multiplying both sides of this inequality by $R^{d\left(1-\frac{1}{p}\right)}$ we get with $C=\omega_d^{1-\frac{1}{p}}$
\bea
C\|u\|_p&=&\left(\int_{\RR^d}u^p(x)\dx\right)^{\frac{1}{p}}\left(\int_{B_R(0)}\dx
\right)^{1-\frac{1}{p}} R^{d\left(\frac{1}{p}-1\right)}\nonumber\\
 &\ge& R^{d\left(\frac{1}{p}-1\right)}\int_{B_R(0)}u(x)\dx\ge R^{d\left(\frac{1}{p}-1\right)}M-R^{-\gamma}w_\gamma.\nonumber
 \eea 
 Taking the optimal $R$, {\it i.e.}  $R^\gamma=Cw_\gamma/M$, we get \rf{m-L}.

Now, it is clear that if condition \rf{cond:blow-up} for blow up is satisfied for some $\gamma\in(1,2]$ and a suitable constant $c$, 
then for $p=d/(\alpha+\beta-2)$ appearing in Theorem \ref{existence}.ii,
we have $\|u\|_p\ge CM M^{-\frac{\gamma}{d+2-\alpha-\beta}\frac{d}{\gamma}\left(1-\frac{1}{p}\right)}=C$, 
so condition \rf{smallness} for global existence is violated for sufficiently small $\varepsilon>0$. 
{\it Vice versa,} if \rf{smallness} is satisfied, \rf{cond:blow-up} cannot be true with small constants $c$. 
}
\end{rem}
\medskip

\begin{rem}{\rm  Note that one can prove in a similar way the inequality 
 \be
\|u\|_{M^p}\ge C M \left(\frac{M}{w_\gamma}\right)^{\frac{d}{\gamma}\left(1-\frac{1}{p}\right)},\label{Morrey}
\ee 
where the Morrey space norm $\|u\|_{M^p}\ (\le C\|u\|_p)$ is defined as 
$$\sup_{R>0,\, x_0\in\RR^d}R^{d\left(\frac{1}{p}-1\right)}\int_{B_R(x_0)}|u(x)|\dx,$$
see \cite[(15)]{B3} in the particular case $\gamma=2$.  
The scale of Morrey spaces is relevant to study mean field type problems related to the Keller--Segel model, since the Morrey norms are particularly well suited to measure the (local) concentration of densities, see {\it e.g.} \cite{B-SM}.}
\end{rem}
\medskip

We prove Theorem \ref{main} and Corollary \ref{cor:main} by showing 
the extinction in a~finite time of the function 
$w(t)=\int_{\RR^d}\varphi_\gamma (x) u(x,t)\dx$ where 
$\varphi_\gamma(x)$ is smooth and behaves like  $|x|^\gamma$ with some 
$\gamma \in (1,\alpha)$ for large $|x|$, see \rf{phi}  below. 
Note that if $\alpha<2$, we cannot expect the existence of higher order moments 
$w_\gamma$ defined in \rf{m-gamma} 
with $\gamma\ge\alpha$. Indeed, even for the linear equation $\partial_tv+(-
\Delta)^{\alpha/2}v=0$, the fundamental 
solution $p_\alpha(x,t)$ behaves like $p_\alpha(x,t)\sim 
\left(t^{d/\alpha}+|x|^{d+\alpha}/t\right)^{-1}$, and therefore the moment 
\rf{m-gamma} with $\gamma\ge\alpha$ cannot be finite, see 
\cite{BK08} and references given there.
Thus, we cannot apply the usual reasoning which involves an analysis of the 
evolution of the second moment $w_2$ of the solution  because the integral defining $w_2$ may diverge. 

We recall  in Proposition \ref{prop1} a result showing  that the moment \rf{m-gamma} is finite
for a large class of initial conditions in the case of $\gamma <\alpha$.

\section{Existence of solutions}\label{sec:exist}
\setcounter{equation}{0}

We are going to construct solutions to system \rf{eq}--\rf{BB} via the 
following integral formulation
\begin{equation}\label{duhamel}
u(t)=S_\alpha(t)u_0 -\int_0^t \nabla \cdot S_\alpha(t-\tau)\big(u(\tau)Bu(\tau)\big)
\,{\rm d\tau},
\end{equation}
{\it i.e.}, we  consider {\em mild} solutions. 
Here $S_\alpha(t)u_0=p_\alpha(t)*u_0$ is the solution to the 
 linear Cauchy problem 
\begin{equation}
\label{leq}
 \partial_t v+(-\Delta)^{\alpha/2}v=0, \quad v(x,0)=u_0, 
\end{equation}
and   $p_\alpha(x,t)$
is the fundamental solution of \rf{leq} which can be represented via the Fourier transform
$ \widehat p_\alpha(\xi,t)=e^{-t|\xi|^\alpha}$.
In particular,
$$
p_\alpha(x,t)=t^{-d/\alpha}P_\alpha(xt^{-1/\alpha}),
$$
where $P_\alpha$ is the inverse Fourier transform of $e^{-
|\xi|^\alpha}$, see \cite[Ch.~3]{J} and \cite{BK08} for more details.
It is well known that for every $\alpha\in (0,2)$ the function
$P_\alpha$
is  smooth, nonnegative, and satisfies the (optimal) estimates
\begin{equation}
\label{EP}
0<P_\alpha(x)\le C(1+|x|)^{-(\alpha+d)} \;\; \hbox{and}\;\; 
|\nabla P_\alpha(x)|\le C(1+|x|)^{-(\alpha+d+1)}
\end{equation}
for a constant $C$ and all $x\in\RR^d$.
Hence, it follows immediately from the Young inequality for the convolution
and from the self-similar form of the kernel $p_\alpha(x,t)$ that for every
$1\leq q\leq p\leq \infty$ there exists $C=C(p,q,\alpha)>0$ such that
\begin{equation}\label{Sat:1}
\|S_\alpha(t)u_0\|_p\leq Ct^{-\frac{d}{\alpha}\left(\frac{1}{q}-\frac{1}{p}\right)}\|u_0\|_q
\end{equation}
and
\begin{equation}\label{Sat:2}
\|\nabla S_\alpha(t)u_0\|_p\leq Ct^{-\frac{d}{\alpha}\left(\frac{1}{q}-\frac{1}{p}\right)-\frac{1}{\alpha}}\|u_0\|_q
\end{equation}
for every $u_0\in L^q(\RR^d)$ and all $t>0$.

\bigskip

The construction of solution to \rf{duhamel} in $L^p$ spaces 
 is based on the following abstract result, see {\it e.g.} \cite{LR}, \cite{Mey99}.

\begin{lem}
\label{lem:xyB}
Let $(\X, \|\cdot\|_\X)$ be a Banach space and $H:\X\times
\X\to
\X$ a bounded bilinear form satisfying
$
\|H(x_1,x_2)\|_\X\leq \eta \|x_1\|_\X \|x_2\|_\X
$ for all $x_1,x_2\in\X$ and a constant $\eta>0$.
Then, if $0<\varepsilon<1/(4\eta)$ and if $y\in\X$ is such that
$ \|y\|<\varepsilon$, the equation
$u=y+H(u,u)$ has a solution in $\X$ such that $\|u\|_\X\leq
2\varepsilon$.
This solution is the
only one in the ball $\bar B(0,2\varepsilon)$.
\end{lem}
We skip an  easy proof of this lemma which  is a~direct consequence of the Banach fixed point theorem.

\bigskip 
\noindent {\it Sketch of the proof of Theorem  \ref{existence}.}
The proof consists in constructing solutions
 to the ``quadratic'' equation \rf{duhamel} using Lemma \ref{lem:xyB}
with
$y=S_\alpha(t)u_0$ and with the bilinear form 
\begin{equation}\label{def:H}
H(u,v)= -\int_0^t \nabla \cdot S_\alpha(t-\tau)\big(u(\tau)Bv(\tau)\big)
\,{\rm d\tau}.
\end{equation}

\medskip

{\it Local existence of solutions.} It suffices to obtain estimates  
required by Lemma \ref{lem:xyB} in the Banach space $\X_T^p=C([0,T], L^p(\RR^d))$ supplemented with the usual norm
$\|u\|_{\X_T^p}=\sup_{t\in [0,T]}\|u(t)\|_p$. 

By inequality \rf{Sat:1}, we immediately obtain $\|S_\alpha(\cdot)u_0\|_{\X_T^p}\leq \|u_0\|_p$. 

Combining  the definition \rf{BB} of the operator $B$ with the Hardy--Little\-wood--Sobolev inequality we obtain
\begin{equation}\label{HSL}
\|B(u)\|_q\leq s_{d,\beta} \big\||\cdot|^{-d+\beta-1}*u\big\|_q\leq C\|u\|_p
\end{equation}
for every $1<p<q<\infty$ satisfying 
$\frac{1}{p}-\frac{\beta-1}{d}=\frac{1}{q}$. Moreover, inequality
\rf{Sat:2} and the H\"older inequality lead to
\begin{equation}\label{Sat:2+Holder}
\begin{split}
\|\nabla S_\alpha(t-\tau)(uB(v))\|_p&\leq C(t-\tau)^{-\frac1\alpha-\frac d\alpha
\left(\frac1r-\frac1p\right)}\|uB(v)\|_r\\
&\leq C(t-\tau)^{-\frac1\alpha-\frac d\alpha
\left(\frac1r-\frac1p\right)}\|u\|_p\|B(v)\|_q.
\end{split}
\end{equation}
where 
$\frac{1}{r}=\frac{1}{p}+\frac{1}{q}$. Hence, by inequalities \rf{HSL} and \rf{Sat:2+Holder}, there exists a~constant $C$ such that for all $u,v \in \X^p_T$, 
the bilinear form
 \rf{def:H} satisfies
\begin{equation}\label{H:est}
\begin{split}
\|H(u,v)\|_{\X_T^p}&\leq C\sup_{t\in [0,T]}\int_0^t 
(t-\tau)^{-\frac1\alpha-\frac d\alpha
\left(\frac1r-\frac1p\right)}\|u(\tau)\|_p\|v(\tau)\|_p\,{\rm d\tau}\\
&\leq CT^{1-\frac1\alpha-\frac d\alpha
\left(\frac1r-\frac1p\right)} \|u\|_{\X_T^p}\|v\|_{\X_T^p}.
\end{split}
\end{equation}
In estimates \rf{H:est}, we have used the relations
$
\frac1r=\frac1p+\frac1q= \frac2p-\frac{\beta-1}d,
$
and we assume that 
\begin{itemize}
\item $p>d/(\alpha+\beta-2)$ in order to have $1/\alpha+(d/\alpha)(1/r-1/p)<1$;

\item $p>2d/(d+\beta-1)$ to guarantee that $r>1$;

\item $p\leq d/(\beta-1)$ to be sure that $r\leq p$.

\end{itemize}

Choosing a sufficiently small $T>0$ in \rf{H:est} we complete the proof 
of Theorem \ref{existence}.i by an application of Lemma \ref{lem:xyB}.

\medskip

{\it Global in time  solutions.}
Here, the reasoning is completely analogous: using inequalities
\rf{Sat:1}, \rf{Sat:2}, \rf{HSL}--\rf{H:est}
we estimate the bilinear form
\rf{def:H} in the Banach space
\begin{equation*}
\begin{split}
\X^p=&C([0,\infty), L^{d/(\alpha+\beta-2)}(\RR^d))\\
&\cap \{u\in C((0,\infty),L^p(\RR^d))\;:\;
\sup_{t>0} t^{\frac{d}{\alpha}\left(\frac1p-\frac{\alpha+\beta-2}d\right)} \|u(t)\|_p<\infty\}
\end{split}
\end{equation*}
 supplemented with the norm
$$
\|u\|_{\X^p}=\sup_{t>0}\|u(t)\|_{d/(\alpha+\beta-2)} +
\sup_{t>0} t^{\frac{d}{\alpha}\left(\frac1p-\frac{\alpha+\beta-2}d\right)} \|u(t)\|_p.
$$
We skip further details of this standard reasoning.

\medskip
{\it Nonnegativity property.}
In order to prove  that $u_0\geq 0$ implies $u(t)\geq 0$, it suffices
 to study the function 
$u^-(x,t)= \max \{-u(x,t), 0\}$ and to follow the  arguments
either  from \cite[Prop.~3.1]{BKW2} or from \cite[Prop.~2]{DI} in order
to show that $u^-(x,t)\equiv 0$. 
Here, we do not give a detailed presentation  because the 
proof from \cite[Lemma 2.7]{LRZ08} can be rewritten in this more general case.

\medskip

{\it Conservation of the integral.}
If $u_0\in L^1(\RR^d)$, one should repeat the fixed point argument from i) (or
from ii))
in the space $\X_T^p\cap C([0,T], L^1(\RR^d))$ 
(in $\X^p\cap C([0,\infty), L^1(\RR^d))$, resp.)  
in order to have $u(t)\in L^1(\RR^d)$ for all $t\in [0,T]$ ($t\in (0,\infty)$, resp.). Next, it suffices to integrate over $\RR^d$ the both sides of equation \rf{duhamel}. Using the following consequences of the Fubini theorem
 $$
\int_{\RR^d} p_\alpha(t)*u_0(x)\dx=\int_{\RR^d} u_0(x)\dx
$$
and
$$
\int_{\RR^d} \nabla p_\alpha(t)*v(x)\dx=0\quad \mbox{for every}\quad
v\in L^1(\RR^d),
$$
we conclude that $\int_{\RR^d} u(x,t)\dx =\int_{\RR^d} u_0(x)\dx$ for every
$t\in [0,T]$.
\qed

\begin{rem}{\rm
The reasoning from the proof of Theorem \ref{existence}  
follows the lines of 
the usual proof of the local in time existence (as well as
the global in time existence for small initial conditions)
 of solutions to system \rf{eq}--\rf{BB} with $\alpha=2$. 
Since that argument is based on an integral representation 
analogous to that in \rf{duhamel} and on counterparts 
 of decay estimates from \rf{Sat:1}--\rf{Sat:2}, it can be easily 
adopted to the more general case of $\alpha\in (1,2]$.
Examples of such a reasoning applied to various semilinear models and 
 realized in  miscellaneous Banach spaces 
can be found in \cite{B-SM,BB08, BCGK, BWu,K99, KS08, LR, Mey99}.}
\end{rem}

\bigskip

Next, we recall weighted  estimates of solutions
to the linear Cauchy problem \rf{leq} which have been proved in, {\it e.g.},  
\cite{BK08}.
In the following, we  use the weighted $L^\infty$ spaces 
\begin{equation*}\label{La}
L^\infty_{\vartheta}(\RR^d)= \{v\in L^\infty(\RR^d)\,:\,
\|v\|_{L^\infty_\vartheta}\equiv
{\rm ess\, sup}_{x\in\RR^d}(1+|x|)^{\vartheta}|v(x)|<\infty\}
\end{equation*}
for a fixed $\vartheta\ge0$.

\begin{lem}{\rm (\cite[Lemma 3.1]{BK08})}
\label{LL1}
Assume that $v_0\in L^\infty_{\alpha+d}(\RR^d)$.
There exists $C>0$ independent of $v_0$ and $t$ such that
\begin{eqnarray*}
\label{St 2}
\|S_\alpha(t)v_0\|_{L^\infty_{\alpha+d}}&\le&
C(1+t)\|v_0\|_{L^\infty_{\alpha+d}} \,,\\
\label{St 3}
\|\nabla S_\alpha(t) v_0\|_{L^\infty_{\alpha+d}} &\le&
Ct^{-1/\alpha}\|v_0\|_{L^\infty_{\alpha+d}}
    +Ct^{1-1/\alpha}\|v_0\|_1 \,,
\end{eqnarray*}
\end{lem}
\bigskip

With these estimates, we are in a position to construct solutions to system
\rf{eq}--\rf{BB} also in the weighted space 
$L^\infty_{\alpha+d}(\RR^d)\subset L^1(\RR^d)\cap L^\infty(\RR^d)$.

\begin{prop}
\label{prop1}
Let $\alpha\in(1,2]$ and $b\in (1,d]$.
Assume that $u$ is the solution of system 
\eqref{eq}--\eqref{BB}, constructed in Theorem \ref{existence},
supplemented with the initial condition 
 $u_0\in L^\infty_{\alpha+d}(\RR^d)$.  
Then
$
u\in C([0,T], L^\infty_{\alpha+d}(\RR^d)).  
$
\ 
In particular, for each $\gamma<\alpha$ we have $\int_{\RR^d}|x|^\gamma u(x,t)\dx<\infty$. 
\end{prop}

\noindent {\it Sketch of proof.}
Here, it suffices to modify slightly the argument from \cite[Proof of 
Prop. 3.3.i]{BK08} written the  case of a convection equation
with the fractional Laplacian. 
In that reasoning (as well as in  the above proof of Theorem \ref{existence}),
we construct solutions to equation \rf{duhamel} by applying Lemma
\ref{lem:xyB} with the Banach space 
$\X_T=C([0,T], L^\infty_{\alpha+d}(\RR^d))$. 
Obviously, $y=S_\alpha(\cdot)u_0\in \X_T$ by Lemma \ref{LL1}. 
Next, we show the following estimate of the bilinear form from \rf{def:H}
$$
\|H(u,v)\|_{\X_T}\leq C T^{1-1/\alpha}\|u|_{\X_T}\|v\|_{\X_T}
$$ 
for all $u,v\in \X_T$ and a constant $C$ independent of $u,v$.
Here, in order to estimate  the singular integral operator \rf{BB}
in the spaces $L^\infty_{\vartheta}(\RR^d))$,  
one should follow the reasoning from \cite[Sec. 2]{BB08}. 
Let us omit other details of this classical argument.
\qed

\section{Blow up of solutions}
\setcounter{equation}{0}

The main role in our proof of the blow up of solutions to \rf{eq}--\rf{BB}
is played by the 
following smooth nonnegative weight function  on $\RR^d$
\begin{equation}\label{phi}
\f(x)=\f_\gamma(x)\equiv(1+|x|^2)^{\gamma/2}-1 
\end{equation}
with $\gamma\in(1,2]$. 
Since $(1+|x|^2)^\gamma\le (1+|x|^\gamma)^2$, we have for each $\varepsilon>0$, suitably chosen $C(\varepsilon)>0$,
 and for every $x\in \RR^d$ 
\begin{equation}\label{phi:est}
\f(x)\leq |x|^\gamma\leq \varepsilon +C(\varepsilon)\f(x).
\end{equation}

Next, let us  state two auxiliary results concerning the weight function $\f$ 
which will be used in the proof
 of Theorem \ref{main}. Here, for a given $\f\in C^2(\RR^d)$, we denote by 
$D^2\f$ its Hessian matrix. Moreover, the scalar product of vectors 
$x,y\in \RR^d$ is denoted by $x\cdot y$. If $A$ is either a vector or a matrix,
the expression $|A|$ means its Euclidean norm.

\begin{lem}\label{finite}
Let $\alpha\in (1,2)$,  
$\gamma \in (1,\alpha)$, and $\f$ be defined by \rf{phi}.
Then 
\be
(-\Delta)^{\alpha/2}\f\in L^\infty(\RR^d).\label{fin}
\ee
\end{lem}

\proof
First note that  by a direct computation we have
\begin{equation}
\nabla \f(x)= 
\gamma(1+|x|^2)^{\frac{\gamma}{2}-1}x\label{p1}
\end{equation}
and
\begin{equation}\label{p2}
\partial_{x_j}\partial_{x_i} \f(x)=  
\left(\gamma(1+|x|^2)\delta_{i\, j}-\gamma(2-\gamma)x_i\,x_j\right)(1+|x|^2)^{\frac{\gamma}{2}-2}. 
\end{equation}
In particular,  for every $R>0$ there exists $C(R,\gamma)>0$ 
such that for all
$|x|\geq R$ we have
\begin{equation}\label{der:est}
|\nabla \f(x)|\leq C(R,\gamma)|x|^{\gamma -1} \quad\mbox{and}\quad
|D^2\f(x)|\leq C(R,\gamma)|x|^{\gamma-2}.
\end{equation}

Now, we apply  the  following L\'evy--Khintchine
integral representation of the fractional Laplacian
\begin{equation}
(-\Delta)^{\alpha/2}\f(x)=C(d,\alpha)
\int_{\RR^d}\frac{\f(x+y)-\f(x)-\nabla\f(x)\cdot y}{|y|^{d+\alpha}}\dy
\label{repr}
\end{equation}
(with a suitable constant $C(d,\alpha)$) which is valid for every 
$\alpha\in (1,2)$, 
see, {\it e.g.}~\cite[Th.~1]{DI} for a detailed proof of that version of L\'evy--Khintchine formula.
Using the Taylor expansion and estimates \rf{der:est} one can immediately show
that $(-\Delta)^{\alpha/2}\f(x)$ is well defined for every $x\in\RR^d$ and, 
moreover, $\sup_{|x|\leq R}|(-\Delta)^{\alpha/2}\f(x)|<\infty$ for each $R>0$.

In order to obtain an estimate uniform in $x\in \RR^d$, we assume that 
$|x|\geq 1$ and we shall estimate the integral on the right hand side of 
\rf{repr} for $|y|\leq |x|/2$ and $|y|>|x|/2$, separately.

If $|y|\leq |x|/2$, by the Taylor formula and the second inequality in 
\rf{der:est}, we obtain
\begin{equation*} 
\begin{split}
|\f(x+y)-\f(x)-\nabla\f(x)\cdot y|&\leq \frac12|y|^2 \int_0^1 
|D^2 \f(x+sy)|\,{\rm ds}\\ 
&\leq C|y|^2 \int_0^1 |x+sy|^{\gamma-2}\,{\rm ds}.
\end{split}
\end{equation*}
Since $|y|\leq |x|/2$ and $s\in [0,1]$ we can estimate
$$
|x+sy|\geq \big||x|-s|y|\big|\geq \big||x|-|y|\big|\geq \frac12|x|.
$$
Consequently, for $\gamma-2<0$, we obtain 
\begin{equation}\label{ineq:1}
\begin{split}
\Bigg|\int_{|y|\leq |x|/2}&\frac{\f(x+y)-\f(x)-\nabla\f(x)\cdot y}{|y|^{d+\alpha}}
\dy\Bigg|\\
&\hspace{2cm}\leq
C|x|^{\gamma-2}\int_{|y|\leq |x|/2} \frac{\dy}{|y|^{d+\alpha-2}}
=C|x|^{\gamma-\alpha}
\end{split}
\end{equation}
for all $|x|\geq 1$ and a constant $C>0$ independent of $x$.

If $|y|\geq |x|/2$ and $|x|\geq 1$, we combine first inequality 
from \rf{der:est} (remember that $\gamma-1>0$) with the Taylor expansion
to show 
\begin{equation*}
|\f(x+y)-\f(x)|\leq |y|\int_0^1|\nabla\f(x+sy)|\,{\rm ds}
\leq C|y|\left(|x|^{\gamma-1}+|y|^{\gamma-1}\right).
\end{equation*}
Hence, 
\begin{equation}\label{ineq:2}
\begin{split}
\Bigg|&\int_{|y|> |x|/2}\frac{\f(x+y)-\f(x)-\nabla\f(x)\cdot y}{|y|^{d+\alpha}}
\dy\Bigg|\\
&\leq
C \Bigg( |x|^{\gamma-1}\int_{|y|> |x|/2} \frac{\dy}{|y|^{d+\alpha-1}}
+\int_{|y|> |x|/2} \frac{\dy}{|y|^{d+\alpha-\gamma}}\Bigg)
=C|x|^{\gamma-\alpha}
\end{split}
\end{equation}
for all $|x|\geq 1$ and a constant $C>0$ independent of $x$. 

Finally, inequalities \rf{ineq:1} and \rf{ineq:2} complete the proof because 
$\gamma<\alpha$.
\qed

\begin{rem} 
{\rm 
Note that above 
we have, in fact,  proved that  
$$
\sup_{x\in\RR^d}\left(1+|x|^{\alpha-\gamma}\right)\left|(-
\Delta)^{\alpha/2}\f(x)\right|<\infty\quad
\mbox{for every}\ \gamma\in (1,\alpha).
$$
} 
\end{rem}

\bigskip

\begin{lem}\label{weight}
 For every $\gamma\in (1,2]$, the function 
$\f$  defined in \rf{phi} is  locally uniformly convex  on $\RR^d$. Moreover, there 
exists $K=K(\gamma)$ such that the following inequality
\be
\left(\nabla\f(x)-\nabla\f(y)\right)\cdot (x-y)\ge \frac{K|x-y|^2}{1+|x|^{2-
\gamma}+|y|^{2-\gamma}} \label{grad}
\ee
holds true for  all $x,\, y\in\RR^d$. 
\end{lem}

\proof
Using the explicit expression for the Hessian matrix of $\f$ in \rf{p2} we obtain
 \begin{equation}\label{Jacobian}
D^2\f(x)y\cdot y=
\frac{ \gamma(1+|x|^2)|y|^2 -\gamma(2-\gamma)
\left(\sum_ix^2_iy^2_i +\sum_{i\neq j} x_ix_jy_iy_j\right)}{(1+|x|^2)^{2-\gamma/2}}
\end{equation}
for every $x,y\in \RR^d$. 
Now, by the elementary inequality 
$x_ix_jy_iy_j\leq \frac12(x_i^2y_j^2+x_i^2y_j^2)$ we immediately obtain
$$
\sum_{i\neq j}x_ix_jy_iy_j \leq \sum_{i\neq j} x_i^2y_j^2.
$$
Consequently,
\begin{equation}\label{sum}
\sum_{i}x_i^2y_i^2 +\sum_{i\neq j} x_ix_jy_iy_j\leq \sum_{i,j}x_i^2y_j^2=|x|^2|y|^2.
\end{equation}
Since $\gamma \in (1,2]$, applying estimate \rf{sum} to \rf{Jacobian} we 
get the inequality
$$
D^2\f(x)y\cdot y\geq \frac{\gamma(1 +(\gamma-1)|x|^2)|y|^2}{(1+|x|^2)^{2-\gamma/2}}
$$
which leads directly to the estimate from below 
\begin{equation}\label{Jacobian2}
D^2\f(x)y\cdot y\geq \frac{(\gamma-1)|y|^2}{(1+|x|^2)^{1-\gamma/2}}.
\end{equation} 
Finally, it follows from the integration of the second derivative of $\f$ that 
\begin{equation*}
\Big(\nabla\f(x)
-\nabla\f(y)\Big)\cdot (x-y)
=\int_0^1 D^2\f\big(x+s(y-x)\big)(x-y)\cdot (x-y) \ds. 
\end{equation*}
Hence, using inequality \rf{Jacobian2}
and the estimate
$$
(1+|x+s(y-x)|^2)^{1-\gamma/2}\leq C(1+|x|^{2-\gamma}+|y|^{2-\gamma}),
$$
valid for all $x,y\in \RR^d$, $s\in [0,1]$ and a constant $C>0$ independent
of $x,y,s$, one can easily complete the proof
of Lemma~\ref{weight}.
\qed
\bigskip


\noindent {\it Proof of Theorem  \ref{main}.}
We consider the function 
$$w=w(t)\equiv\int_{\RR^d}\f(x)\,u(x,t)\dx,$$ 
where $\f$ is defined in \rf{phi} and 
 $1<\gamma<\alpha$.
Note that, in view of inequalities 
\rf{phi:est}, the quantity 
$w$  is essentially equivalent to the  moment $w_\gamma$ of order $\gamma$ of the solution $u$.
 Moreover, it satisfies the 
relation 
\begin{equation}\label{ww}
\begin{split}
\frac{{\rm d}}{{\rm d}t}w=&-\int_{\RR^d}(-\Delta)^{\alpha/2}u(x,t)\,\f(x)\dx -
\int_{\RR^d}u(x,t)Bu(x,t)\cdot\nabla\f(x)\dx\\
=&-\int_{\RR^d}(-\Delta)^{\alpha/2}\f(x)\,u(x,t)\dx\\ 
&- \frac{s_{d,\beta}}2\int_{\RR^d}\int_{\RR^d}\Big(\nabla\f(x)-\nabla\f(y)\Big)\cdot (x-y) 
\frac{u(x,t)u(y,t)}{|x-y|^{d-\beta+2}}\dx \dy
\end{split}
\end{equation}
after using the definition of the form $B$ in \rf{BB} and the symmetrization of 
the double integral. This computation resembles the usual proof of blow up 
involving the second  moments, cf. \cite{B3,BW, CPZ,CPS}. 
\bigskip

i) For $\alpha=2=\gamma$ (hence for $\f(x)=|x|^2$) 
and for $\beta=d$, the equality  \rf{ww} can be rewritten as follows 
$$\frac{{\rm d}}{{\rm d}t}w(t)= 2d M-s_{d,\beta} M^2.$$ 
Evidently, for $M>2d/s_{d,\beta}$, this implies the 
equality $w(T)=0$ for some $T>0$, a contradiction with the global 
existence of nonnegative solutions. Thus, we recover the result in \cite[Prop. 
4.1]{BW}  refined in \cite{CPZ, CPS}.  
\bigskip

ii) For $1<\beta\le d$ and fixed $M>0$, we are going to use the following simple 
identity
\begin{equation*}
\begin{split}
M^2&=\int_{\RR^d}\int_{\RR^d}u(x,t)u(y,t) \dx \dy\nonumber\\ 
&=\int_{\RR^d}\int_{\RR^d}u(x,t)u(y,t)\frac{|x-y|^\nu}{\left(1+|x|^{2-
\gamma}+|y|^{2-\gamma}\right)^{\delta}}\frac{\left(1+|x|^{2-\gamma}+|y|^{2-
\gamma}\right)^{\delta}}{|x-y|^{\nu}}\dx \dy\nonumber
\end{split}
\end{equation*}
with some $\nu>0$ and $\delta>0$. 
We apply now the H\"older inequality with the  
powers $p>1$ and $p'=\frac{p}{p-1}$ chosen so that 
\begin{equation}\label{exp}
\nu p=d-\beta,\ \ \ \delta p=1,\ \ \ \nu p'+(2-\gamma)\delta{p'}=\gamma.
\end{equation} 
Of course, such a choice of $\nu$, $\delta$, $p$ is possible whenever $\beta<d$ 
and $\gamma<2$ because we only need $d-\beta+2-\gamma=\gamma(p-1)$. If 
$\beta=d$, it suffices to take $\nu=0$ and $p=2/\gamma>1$. 
As a consequence, we get
\begin{equation}
\begin{split}
M^2& \leq J(t)^{1/p}\label{M:est1:1}\\
& \times\left(\int_{\RR^d}\int_{\RR^d} u(x,t)u(y,t)|x-y|^{\nu 
p'}\left(1+|x|^{2-\gamma} +|y|^{2-\gamma}\right)^{\delta p'}\dx 
\dy\right)^{1/p'},
\end{split}
\end{equation}
where the integral $J(t)$ satisfies
\begin{equation}\label{J:est}
\begin{split}
J(t)&=\int_{\RR^d}\int_{\RR^d}\frac{u(x,t)u(y,t)}{|x-y|^{d-\beta}} 
\frac{\dx \dy}{1+|x|^{2-\gamma}+|y|^{2-\gamma}}\\
&\leq \frac{1}{K}
\int_{\RR^d}\int_{\RR^d}\Big(\nabla\f(x)-\nabla\f(y)\Big)\cdot (x-y) 
\frac{u(x,t)u(y,t)}{|x-y|^{d-\beta+2}}\dx \dy
\end{split}
\end{equation}
by Lemma \ref{weight}.

It follows from relations \rf{exp} and 
inequalities \rf{phi:est} that there exists a~constant  
$C_1>0$ such that 
\be
|x-y|^{\nu 
p'}\left(1+|x|^{2-\gamma} +|y|^{2-\gamma}\right)^{\delta p'}
\leq C_1\Big(1+\f(x)+\f(y)\Big).\label{phi-phi}
\ee
Hence, \rf{M:est1:1} implies 
\begin{equation}\label{M:est:2}
M^2\leq C_1^{1/p'} 
J(t)^{1/p}\Big( M^2+2Mw(t))\Big)^{1/p'}.
\end{equation}
Going back to identity \rf{ww} we obtain from Lemma \ref{finite} and
from inequalities \rf{J:est}--\rf{M:est:2} that
\begin{equation}
\frac{{\rm d}}{{\rm d}t}w(t) \leq C_2M-
C_3\frac{M^{2p}}{( M^2+2Mw(t))^{p/p'}}
\label{evo}
\end{equation}
with  $C_2=\|(-\Delta)^{\alpha/2}\f\|_\infty$ and a suitable constant $C_3>0$. 

Now, we fix  for a while  $M=M_0$ in \rf{evo} so large in order to have 
\begin{equation}\label{M:est}
C_2M_0-
C_3\frac{M_0^{2p}}{( M_0^2)^{p/p'}}<0.
\end{equation}
Hence, there exists 
$C_4=C_4(M_0)>0$ such that for $0<w(0)\leq C_4$ we still have 
$$
C_2M_0-
C_3\frac{M_0^{2p}}{( M_0^2+2M_0w(0))^{p/p'}}<0.
$$
It is clear that if initially $0<w(0)\leq C_4$ then,
by inequality \rf{evo} with $M=M_0$, the function 
$w(t)$  is decreasing in time. Moreover, 
$$
\frac{{\rm d}}{{\rm d}t}w(t)\le 
C_2M_0-
C_3\frac{M_0^{2p}}{( M_0^2+2M_0w(0))^{p/p'}}<0
$$ 
and, consequently, $w(T)=0$ for some $0<T<\infty$. 
This contradicts the global in 
time existence of regular nonnegative solutions of \rf{eq}--\rf{BB}. 
Finally, note that due to the first inequality in \rf{phi:est}, it suffices to assume
\begin{equation}\label{u0:blow:1}
w(0)\leq \int_{\RR^d} |x|^\gamma u_0(x) \dx\leq C_4 \quad \mbox{and} \quad 
\int_{\RR^d}u_0(x)\dx=M_0,
\end{equation}
in order to obtain the blow up in a finite time of the corresponding 
solution.

Now, assume that 
$\int_{\RR^d}u(x,t)\dx=\int_{\RR^d}u_0(x)\dx=M\neq M_0$.
Recall that system \rf{eq}--\rf{BB} is invariant under the rescaling \rf{scal}.
Choosing $\lambda^{\alpha+\beta-2-d}=M_0/M$ we obtain
$\int_{\RR^d}u^\lambda(x,t)\dx=\int_{\RR^d}u_0^\lambda(x)\dx= M_0$
and, by \rf{u0:blow:1},  the blow up of the solution takes place under the assumption
$$
\int_{\RR^d}|x|^\gamma u_0^\lambda(x) \dx\leq C_4.
$$
Changing the variables and using the explicit form of $\lambda$ we obtain the
blow up of solutions to \rf{eq}--\rf{BB} under the following assumption on the initial condition 
$$
\int_{\RR^d} |x|^\gamma u_0(x) \dx\leq C_4M_0^{-1+\frac{\gamma}{\alpha+\beta-2-d}}
\Bigg(\int_{\RR^d}u_0(x)\dx
\Bigg)^{1+\frac{\gamma}{d+2-\alpha-\beta}}.
$$
\qed

\bigskip

\noindent {\it Proof of Corollary  \ref{cor:main}.}
We follow the proof of Theorem \ref{main}. In particular, we choose $M_\gamma$ so large
that the inequality \rf{M:est} holds true for all $M_0>M_\gamma$. This leads to the blow up
of the corresponding solution under the assumption \rf{u0:blow:1} imposed on the initial data. 

To complete the proof, we use the scaling argument  again. By \rf{scal}, 
$u^\lambda(x,t)=\lambda^2u(\lambda x,\lambda^2 t)$ is a solution for every $\lambda>0$.
Note now that 
$$
\int_{\RR^2}u^\lambda_0(x)\dx= \int_{\RR^2}u_0(x)\dx\;\mbox{and} \; 
\int_{\RR^2}|x|^\gamma u^\lambda_0(x)\dx= \lambda^{-\gamma}\int_{\RR^2}|x|^\gamma u_0(x)\dx.$$
Hence, each initial data $u_0\in L^1(\RR^2, (1+|x|^\gamma)\dx)$ 
satisfying $\int_{\RR^2}u_0(x)\dx=M_0>M^\gamma$ 
leads to the blow up in a finite time of the corresponding solution because 
 the moment condition in \rf{u0:blow:1} can be satisfied  replacing $u$ by $u^\lambda$  and choosing $\lambda$ large enough. 
\qed

\bigskip

\noindent {\bf Acknowledgements.}
The preparation of this paper was partially supported by the Polish Ministry of 
Science grant N201 022 32/0902, the POLONIUM project \'EGIDE no.~13886SG,  
and by the European Commission Marie Curie Host Fellowship 
for the Transfer of Knowledge ``Harmonic Analysis, Nonlinear
Analysis and Probability''  MTKD-CT-2004-013389. 
The authors are greatly indebted to Tomasz Cie\'slak for pointing them out the preprint
\cite{LRZ08}.



\begin{thebibliography}{99}




\bibitem{B-AMSA} 
{\sc P. Biler}, 
{\sl Local and global solvability of parabolic systems modelling chemotaxis,} 
Adv. Math. Sci. Appl. {\bf 8} (1998), 715--743. 



\bibitem{B3} 
{\sc P. Biler},
{\sl Existence and nonexistence of solutions for a~model of gravitational 
interaction of particles.\ III}, 
Colloq. Math. {\bf 68} (1995), 229--239. 

\bibitem{B-SM} 
{\sc P. Biler}, 
{\sl  The Cauchy problem and self-similar solutions for a nonlinear parabolic equation,} Studia Math. {\bf 114} (1995), 181--205.


\bibitem{B-BCP} 
{\sc P. Biler}, 
{\sl Radially symmetric solutions of a chemotaxis model in the plane -- the supercritical case}, 31--42, in: {\it Parabolic and Navier-Stokes Equations}, 
Banach Center Publications {\bf 81}, Polish Acad. Sci.,  Warsaw, 2008.

\bibitem{BB08}
{\sc P. Biler, L. Brandolese}, 
{\sl On the parabolic-elliptic limit
of the doubly parabolic Keller--Segel system
modelling chemotaxis}, (2008), 1--21, 	{\tt arXiv:0804.1000v1 [math.AP]}.


\bibitem{BCGK}
{\sc P. Biler, M. Cannone, I. Guerra, G. Karch,}
{\sl Global regular and singular solutions for a model of gravitating 
particles}, Mathematische Annalen {\bf 330}  (2004), 693--708. 

\bibitem{BFW} 
{\sc P. Biler, T. Funaki, W. A. Woyczy\'nski}, 
{\sl  Interacting particle approximation for nonlocal quadratic evolution 
problems}, Probab. Math. Stat. {\bf 19} (1999), 267--286. 


\bibitem{BKW1}
{\sc P. Biler, G. Karch, W. A. Woyczy\'nski},
{\sl Critical nonlinearity exponent and self-similar asymptotics for L\'evy 
conservation laws}, 
Ann. Inst. H. Poincar\'e -- Analyse non Lin\'eaire {\bf 18} (2001), 613--637.

\bibitem{BKW2} 
{\sc P. Biler, G. Karch, W. A. Woyczy\'nski}, 
{\sl Asymptotics for conservation laws involving L\'evy diffusion generators},  
Studia Math. {\bf 148} (2001), 171--192. 

\bibitem{BW} 
{\sc P. Biler, W. A. Woyczy\'nski}, 
{\sl Global and exploding solutions for nonlocal quadratic evolution problems}, 
SIAM J. Appl. Math. {\bf 59} (1998), 845--869.

\bibitem{BWu} 
{\sc P. Biler, G. Wu}, 
{\sl Two-dimensional chemotaxis models with fractional diffusion}, 
Math. Methods Appl. Sciences {\bf 32} (2009), 112--126; 
\newline DOI: 10.1002/mma.1036, 2008. 

\bibitem{BDP}
{\sc A. Blanchet, J. Dolbeault, B. Perthame}, 
{\sl Two dimensional Keller--Segel model: Optimal critical mass and qualitative 
properties of the solutions}, 
Electron. J. Diff. Eqns. 2006, {\bf 44}, 1--33. 

\bibitem{BK08}
{\sc L. Brandolese, G. Karch},
{\sl Far field asymptotics of solutions to convection equation with anomalous diffusion}, J. Evolution Equations
{\bf 8} (2008), 307--326.



\bibitem{CPS}
{\sc V. Calvez, B. Perthame, M. Sharifi tabar}, 
{\sl Modified Keller--Segel system and critical mass for the $\log$ 
interaction kernel}, Contemporary Math. {\bf 429}, 45--62 (2007), Stochastic 
analysis and pde; Chen, Gui-Qiang (ed.) et al., AMS. 

\bibitem{CPZ}
{\sc L. Corrias, B. Perthame, H. Zaag}, {\sl Global solutions of some chemotaxis and angiogenesis systems in high space dimensions}, Milan J.  Math. {\bf 72} (2004), 1--29.


\bibitem{DI} {\sc J. Droniou, C. Imbert}, 
{\sl Fractal first order partial differential equations}, 
Arch. Rat. Mech. Anal. {\bf 182} (2006), 299--331.


\bibitem{E}
{\sc C. Escudero}, 
{\sl The fractional Keller--Segel model}, Nonlinearity {\bf 19} (2006), 
2909--2918. 


\bibitem{J}
{\sc N. Jacob}, 
{\sl Pseudo-differential Operators and Markov Processes}, vol. 1: Fourier 
analysis and semigroups, Imperial College Press, London, 2001. 

\bibitem{JL}
{\sc W. J\"ager, S. Luckhaus}, 
{\sl   On explosions of solutions to a system of partial differential equations 
modelling chemotaxis}, 
Trans. Amer. Math. Soc.  {\bf  329} (1992), 819--824. 

\bibitem{K99}
{\sc G. Karch},
{\sl Scaling in nolinear parabolic equations},
J. Math. Anal. Appl.,
{\bf 234} (1999), 534--558.

\bibitem{KS08}
{\sc H. Kozono, Y. Sugiyama}, 
{\sl Local existence and finite time blow-up of solutions in the 2-D 
Keller-Segel system,} J. Evol. Equ. {\bf 8} (2008), 353--378. 

\bibitem{LR} {\sc P.-G. Lemari\'e-Rieusset}, {\sl Recent Development in the Navier--Stokes Problem}, Chapman \& Hall/CRC Press, Boca Raton, 2002. 


\bibitem{LRZ08} {\sc D. Li, J. Rodrigo, X. Zhang}, {\sl Exploding solutions for
a nonlocal quadratic evolution problem}, 1--40, preprint on the webpage\newline
{\tt http://www.warwick.ac.uk/staff/J.Rodrigo/research.html}


\bibitem
{Mey99}
{\sc Y. Meyer}, {\sl Wavelets, paraproducts and Navier--Stokes equations},
    Current developments in mathematics, 1996, Internat. Press, Cambridge, MA
    02238-2872 (1999).

\bibitem
{N2000}
{\sc T. Nagai}, {\sl Behavior of solutions to a parabolic-elliptic system 
modelling chemotaxis}, 
  J. Korean Math. Soc. {\bf 37} (2000), 721--732.
 

\end{thebibliography}
\end{document}